\documentclass{amsart}
\usepackage{amssymb}
\usepackage{amsfonts}
\usepackage{amsmath}
\usepackage{graphicx}
\newtheorem{thm}{Theorem}[section]

\newtheorem{proposition}[thm]{Proposition}

\theoremstyle{definition}

\newtheorem{example}[thm]{Example}

\theoremstyle{remark}
\newtheorem{remark}[thm]{Remark}

\numberwithin{equation}{section}
\def\1{{\rm (1)}}
\def\2{{\rm (2)}}
\def\3{{\rm (3)}}
\def\4{{\rm (4)}}
\def\5{{\rm (5)}}
\begin{document}

%%%%%%%%%%%%%%%%%%%%%%%%%%%%%%%%%%%%%%%%%%%%%%%%%%%%%%%%%
%%%%%%%%%%%%%%%%%%%%%%%%%%%%%%%%%%%%%%%%%%%%%%%%%%%%%%%%%
\title[Note on Star Operations over Polynomial Rings]
{Note on Star Operations over Polynomial Rings}

%%%%%%%%%%%%%%%%%%%%%%%%%%%%%%%%%%%%%%%%%%%%%%%%%%%%%%%%%
%%%%%%%%%%%%%%%%%%%%%%%%%%%%%%%%%%%%%%%%%%%%%%%%%%%%%%%%%

\author{A. Mimouni}

\address{Department of Mathematical Sciences,
King Fahd University of Petroleum \& Minerals, P. O. Box 278,
Dhahran 31261, Saudi Arabia}

\email{amimouni@kfupm.edu.sa}

\date{}

%%%%%%%%%%%%%%%%%%%%%%%%%%%%%%%%%%%%%%%%%%%%%%%%%%%%%%%%%
%%%%%%%%%%%%%%%%%%%%%%%%%%%%%%%%%%%%%%%%%%%%%%%%%%%%%%%%%


%%%%%%%%%%%%%%%%%%%%%%%%%%%%%%%%%%%%%%%%%%%%%%%%%%%%%%%%%
%%%%%%%%%%%%%%%%%%%%%%%%%%%%%%%%%%%%%%%%%%%%%%%%%%%%%%%%%
\begin{abstract}
This paper studies the notions of star and semistar operations over
a polynomial ring. It aims at characterizing when every upper to
zero in $R[X]$ is a $*$-maximal ideal and when a $*$-maximal ideal
$Q$ of $R[X]$ is extended from $R$, that is, $Q=(Q\cap R)[X]$ with
$Q\cap R\not =0$, for a given star operation of finite character $*$
on $R[X]$. We also answer negatively some questions raised by
Anderson-Clarke by constructing a Pr\"ufer domain $R$ for which the
$v$-operation is not stable.
\end{abstract}

\maketitle

%%%%%%%%%%%%%%%%%%%%%%%%%%%%%%%%%%%%%%%%%%%%%%%%%%%%%%%%%
%%%%%%%%%%%%%%%%%%%%%%%%%%%%%%%%%%%%%%%%%%%%%%%%%%%%%%%%%
%%%%%%%%%%%%%%%%%%%%%%%%%%%%%%%%%%%%%%%%%%%%%%%%%%%%%%%%%
%%%%%%%%%%%%%%%%%%%%%%%%%%%%%%%%%%%%%%%%%%%%%%%%%%%%%%%%%
\section{Introduction}\label{Int}

Star operations, such as the $t$-closure, the $w$-closure and the
$v$-closure (or divisorial closure) are an essential tool in modern
multiplicative ideal theory for characterizing and investigating
several classes of integral domains.  For an integral domain $R$
with quotient field $K$, if $f\in K[X]$ is an irreducible
polynomial, then we call the prime ideal $P=fK[X]\cap R[X]$ an upper
to zero. Uppers to zero have been used by many authors to
characterize ring-theoretic properties. For example, a domain $R$ is
Pr\"ufer if and only if $P\not\subseteq M[X]$ for each upper to zero
$P$ in $R[X]$ and each maximal ideal $M$ of $R$, \cite[Theorem
19.15]{G1}. A corresponding result exists for Pr\"ufer
$v$-multiplication domains ($PvMD$s). A domain $R$ is a $PvMD$ if
and only if $R$ is an integrally closed $UMT$-domain,
\cite[Proposition 2.6]{HMM}. $UMT$-domains, that is, domains $R$
such that every upper to zero in the polynomial ring $R[X]$ is
$t$-maximal were introduced in \cite{HZ} and studied in
\cite{DHLRZ}, \cite{DHLZ}, \cite{Fa} and investigated in greater
detail in \cite{FGH}. Recently, E. Houston and M. Zafrullah
introduced and studied the notion of $UMV$-domains, that is, domains
$R$ such that every upper to zero in $R[X]$ is a maximal divisorial
ideal (\cite{HZ2}). Also it is well-known that for each nonzero
fractional ideal $I$ of $R$, $(I[X])_{t}=I_{t}[X]$,
$(I[X])_{v}=I_{v}[X]$ and $(I[X])_{w}=I_{w}[X]$ (\cite[Proposition
4.3]{HH} and every $t$-maximal ideal $Q$ of $R[X]$ such that $Q\cap
R=q\not =0$ is of the form $Q=q[X]$ (\cite[Proposition 1.1]{HZ}. The
purpose of this note is to study when uppers to zero are $*$-maximal
ideals and when $Q=(Q\cap R)[X]$ for a $*$-maximal ideal $Q$ of
$R[X]$ with $Q\cap R\not =0$, for a given star operation of finite
character $*$ on $R[X]$. Our main results assert that for a domain
$R$ and a star operation of finite character $*$ on $R[X]$, every
upper to zero is $*$-maximal if and only if $R$ is a $UMT$-domain
and $t$-Max(R[X])=$*$-Max(R[X]) (Theorem~\ref{MR.9}). Also for a
$*$-maximal ideal $Q$ of $R[X]$ with $Q\cap R\not=0$, $Q=(Q\cap
R)[X]$ if and only if $c(Q)^{\bar{*}}\subsetneq R$ if and only if
$c(g)^{\bar{*}}\subsetneq R$ for each $g\in Q$, where $\bar{*}$ is
the star operation of finite character on $R$ induced by $*$
(Proposition~\ref{MR.3}). This characterization leads us to
construct a star operation of finite character $*$ on $R[X]$ such
that $(Q\cap R)[X]\subsetneq Q$ for some $*$-maximal ideal $Q$ of
$R[X]$ with $Q\cap R\not =0$. We close this note by answering
negatively some questions cited in \cite{AC}. Precisely, we
construct a Pr\"ufer domain $R$ such that the $v$-operation on $R$
is not stable.\\

Throughout $R$ is an integral domain with quotient field $K$,
$\bar{F}(R)$ denotes the set of all nonzero $R$-submodules of $K$,
$F(R)$ denotes the set of all nonzero fractional ideals of $R$, i.
e., $E\in F(R)$ if $E\in \bar{F}(R)$ and $dE\subseteq R$ for some
$0\not =d\in R$ and $f(R)$ the set of all nonzero finitely generated
$R$-submodules of $K$. A semistar operation on $R$ is a map $*:
\bar{F}(R)\longrightarrow \bar{F}(R), E\mapsto E^{*}$ satisfying the
following properties for each $E, F\in \bar{F}(R)$ and each $0\not
=a\in K$:\\
$(*_{1})$ $(aE)^{*}=aE^{*}$;\\
$(*_{2})$ $E\subseteq E^{*}$ and if $E\subseteq F$, then
$E^{*}\subseteq F^{*}$;\\
$(*_{3})$ $E^{**}=E^{*}$.\\
In the particular case where $R^{*}=R$, we say that $*|_{F(R)}$ is a
star operation on $R$.\\
A semistar operation $*$ on $R$ is said to be  of finite character
(or of finite type) if $E^{*}=\displaystyle\bigcup\{ F^{*}/ F\in
f(R), F\subseteq E\}$. If $*$ is a semistar operation on $R$, the
map $*_{f}:\bar{F}(R)\longrightarrow \bar{F}(R), E\mapsto
E^{*_{f}}:=\displaystyle\bigcup \{F^{*}/ F\in f(R), F\subseteq E\}$
is a semistar operation of finite character on $R$ called the
semistar operation of finite character associated to $*$.
Obviously, $*$ is of finite character if and only if $*=*_{f}$. A
nonzero ideal $I$ is said to be a $*$-ideal if $I=I^{*}$ and a
$*$-prime ideal is a prime ideal which is a $*$-ideal. It's
well-known that if $*$ is a star operation of finite character and
$I$ is an integral $*$-ideal, then $I$ is contained in a $*$-prime
ideal and a minimal prime ideal over a $*$-ideal is a $*$-prime
ideal. We recall that the $v$-operation is the largest star
operation on $R$ and the $t$-operation is the largest star operation
of finite character on $R$. Finally, we use the notation $*$-Max(R)
to denote the set of all $*$-maximal ideals of $R$, that is,
$*$-prime ideals of $R$ maximal with respect to inclusion, and for a nonzero
polynomial $f\in R[X]$, $c(f)$ is the content of $f$, that is, the
ideal of $R$ generated by all coefficients of $f$.

%%%%%%%%%%%%%%%%%%%%%%%%%%%%%%%%%%%%%%%%%%%%%%%%%%%%%%%%%%%%%%%%%%%%
%%%%%%%%%%%%%%%%%%%%%%%%%%%%%%%%%%%%%%%%%%%%%%%%%%%%%%%%%%%%%%%%%%%%
%%%%%%%%%%%%%%%%%%%%%%%%%%%%%%%%%%%%%%%%%%%%%%%%%%%%%%%%%%%%%%%%%%%%
\section{Main results}\label{MR}

We start with the following proposition showing that every
semistar operation $*$ on the polynomial ring $R[X]$ induces a
semistar operation on $R$. We often refer to it as $\bar{*}$.\\

\begin{proposition}\label{MR.1} Let $R$ be an integral domain,
$X$ an indeterminate over $R$ and $*$ a semistar operation on
$R[X]$. Then the map $\bar{*} : \bar{F}(R)\rightarrow \bar{F}(R),
E\mapsto E^{\bar{*}}=(E[X])^{*}\cap K$ is a semistar operation on
$R$ and $I[X]^{*}=(I^{\bar{*}}[X])^{*}$ for each nonzero fractional
ideal $I$ of $R$. Moreover, if $*$ is of finite character, then so
is $\bar{*}$.
\end{proposition}

\begin{proof} Most of the verification that $\bar{*}$ is a semistar operation is routine.
We consider $(*_{3})$. Thus, let $E\in \bar{F}(R)$. Then $E^{\bar{*}}[X]=(E[X]^{*}\cap
K)[X]=E[X]^{*}[X]\cap K[X]= E[X]^{*}\cap K[X]\subseteq
E[X]^{*}$. Hence $E[X]\subseteq E^{\bar{*}}[X]\subseteq E[X]^{*}$.
So $E[X]^{*}\subseteq E^{\bar{*}}[X]^{*}\subseteq
E[X]^{**}=E[X]^{*}$. Hence $E[X]^{*}=E^{\bar{*}}[X]^{*}$, and
$E^{\bar{*}\bar{*}}=E^{\bar{*}}[X]^{*}\cap K=E[X]^{*}\cap
K=E^{\bar{*}}$. It follows that $\bar{*}$ is a semistar
operation on $R$.\\

Now, Let $I$ be a nonzero fractional ideal of $R$. By the proof of
$(*_{3})$, $I[X]^{*}=(I^{\bar{*}}[X])^{*}$.\\

Assume that $*$ is of finite character. Let $E\in \bar{F}(R)$ and
let $x\in E^{\bar{*}}$. Then $x\in E[X]^{*}\cap K$. Since $*$ is of
finite character, then there exists a finitely generated
$R[X]$-submodule $F$ of $K(X)$ such that $F\subseteq E[X]$ and $x\in
F^{*}$. Set $F=\displaystyle\sum_{i=1}^{i=n}f_{i}R[X]$ and let
$I=\displaystyle\sum_{i=1}^{i=n}c(f_{i})$. Clearly $I$ is f.g.,
$I\subseteq E$ (since each $f_{i}\in F\subseteq E[X]$) and $x\in
F^{*}\subseteq I[X]^{*}$, since $F\subseteq I[X]$. Hence $x\in
I^{\bar{*}}$ and therefore $\bar{*}$ is of finite character.\\

Finally, if $*$ is a star operation on $R[X]$, then
$R^{\bar{*}}=R[X]^{*}\cap K=R[X]\cap K=R$. Hence $\bar{*}|_{F(R)}$
is a star operation on $R$.
\end{proof}

\begin{remark}\label{MR.2} If $*$ denotes the $t$-, respectively
$v$-,  respectively $w$-operation on $R[X]$, then $\bar{*}$ is the
$t$-, respectively $v$-, respectively $w$-operation on $R$. Indeed,
by \cite[Proposition 4.3]{HH},  for each $I\in F(R)$,
$I^{\bar{*}}:=I[X]^{*}\cap K=I^{*}[X]\cap K=I^{*}$.
\end{remark}

The next proposition characterizes $*$-maximal ideals of $R[X]$ that
are extended from $R$ for a given star operation of finite character
on $R[X]$.\\

\begin{proposition}\label{MR.3} Let  $*$ be a star operation of
finite character on $R[X]$ and $Q$ a $*$-maximal ideal of $R[X]$
such that $Q\cap R=q\not =0$. The following statements are
equivalent:\\
$(i)$ $Q=q[X]$;\\
$(ii)$ $c(Q)^{\bar{*}}\subsetneq R$;\\
$(iii)$ $c(g)^{\bar{*}}\subsetneq R$ for each $g\in Q$.
\end{proposition}

\begin{proof} $(i)\Longrightarrow (ii)$. By Proposition~\ref{MR.1},
it is clear that $I[X]^{*}=R[X]$ if and only if $I^{\bar{*}}=R$. Let
$Q$ be a $*$-maximal ideal of $R[X]$. If $Q=q[X]$, then $q$ is
$\bar{*}$-prime and $c(Q)\subseteq q$. Hence
$c(Q)^{\bar{*}}\subsetneq R$.\\
$(ii)\Longrightarrow (i)$. If $c(Q)^{\bar{*}}\subsetneq R$, then
$c(Q)[X]^{*}\subsetneq R[X]$. Since $Q\subseteq c(Q)[X]\subseteq
c(Q)[X]^{*}$ and $Q$ is $*$-maximal, then $Q=c(Q)[X]=c(Q)[X]^{*}$.
Hence $q=Q\cap R=c(Q)[X]\cap R=c(Q)$. So $Q=c(Q)[X]=q[X]$, as
desired.\\
$(ii)\Longrightarrow (iii)$, Trivial.\\
$(iii)\Longrightarrow (ii)$ Suppose that $c(Q)^{\bar{*}}=R$. Since
$\bar{*}$ is of finite character, then there exists a finitely
generated ideal $I\subseteq c(Q)$ such that $I^{\bar{*}}=R$. So
there exist polynomials $f_{1}, \dots, f_{n}\in Q$ such that
$(c(f_{1})+\dots +c(f_{n}))^{\bar{*}}=R$. Set $r_{i}=deg(f_{i})$ and
let $g=f_{1}+X^{r_{1}+1}f_{2}+X^{r_{1}+r_{2}+1}f_{3}+\dots
+X^{r_{1}+\dots+r_{n-1}+1}f_{n}$. Clearly $g\in Q$ and
$c(g)^{\bar{*}}=R$, a contradiction. It follows that
$c(Q)^{\bar{*}}\subsetneq R$.
\end{proof}

The next theorem is an analogue of \cite[Theorem
1.1]{FGH}.

\begin{thm}\label{MR.4} Let $*$ be a star operation of finite character
on $R[X]$. Consider the following assertions:\\
$(i)$ Every upper to zero is $*$-maximal;\\
$(ii)$ For every upper to zero $Q$ in $R[X]$,
$c(Q)^{\bar{*}}=R$;\\
$(iii)$ For every upper to zero $Q$ in $R[X]$, there exists $g\in Q$
such that $c(g)^{\bar{*}}=R$;\\
$(iv)$ $Q\not\subseteq p[X]$ for each upper to zero $Q$ in $R[X]$
and each $\bar{*}$-prime $p$ of $R$;\\
$(v)$ $Q\not\subseteq M[X]$ for each upper to zero $Q$ in $R[X]$
and each $\bar{*}$-maximal ideal $M$ of $R$.\\
Then $(i)\Longrightarrow (ii)\Leftrightarrow (iii)\Leftrightarrow
(iv)\Leftrightarrow (v)$. Moreover, if $*$ satisfies the
conditions in Proposition~\ref{MR.3} for each $*$-maximal ideal
$Q$ of $R[X]$ with $Q\cap R\not =0$, then all 5 conditions are equivalent.
\end{thm}

\begin{proof} $(i)\Longrightarrow (ii)$ Let $Q$ be an upper to zero
in $R[X]$. Since $Q\subsetneq c(Q)[X]$ and $Q$ is $*$-maximal, then
$(c(Q)[X])^{*}=R[X]$. Hence $c(Q)^{\bar{*}}=(c(Q)[X])^{*}\cap
K=R[X]\cap K=R$.\\

$(ii)\Longrightarrow (iii)$ This follows as in the proof of $(iii)\Longrightarrow (ii)$ of
Proposition~\ref{MR.3}.\\

$(iii)\Longrightarrow (iv)$ Let $p$ be a $\bar{*}$-prime ideal of
$R$ and $Q$ an upper to zero in $R[X]$. If $Q\subseteq p[X]$, then
$g\in p[X]$. Hence $R=c(g)^{\bar{*}}\subseteq p^{\bar{*}}=p$,
which is absurd.\\

$(iv)\Longrightarrow (v)$ Trivial.\\

$(v)\Longrightarrow (ii)$ Let $Q$ be an upper to zero in $R[X]$. If
$c(Q)^{\bar{*}}\subsetneq R$, then $c(Q)\subseteq M$ for some
$\bar{*}$-maximal ideal $M$ of $R$. Hence $Q\subseteq
c(Q)[X]\subseteq M[X]$, which is absurd.\\

Now assume that $N=(N\cap R)[X]$ for each $N\in *$-$Max(R[X])$ with
$N\cap R\not =0$. Under this assumption, we show that $(v)\Longrightarrow (i)$.
Let $Q$ be an upper to zero in $R[X]$. Since $Q$
is a $t$-prime ideal, then $Q$ is a $*$-ideal. So $Q\subseteq N$ for
some $*$-maximal ideal $N$ of $R[X]$. If $Q\subsetneq N$, then
$M=N\cap R\not =0$. So $M$ is $\bar{*}$-maximal and $Q\subseteq
N=M[X]$, a contradiction. Hence $Q=N$ is $*$-maximal as desired.
\end{proof}

Below, we give an example showing that a $*$-maximal ideal $Q$ of
$R[X]$ with $Q\cap R\not =0$ need not be extended from $R$ and hence
the statements of Theorem~\ref{MR.4} are not equivalent (unlike the situation in \cite[Theorem
1.1]{FGH}). First we need the following two propositions in which we construct a star
operation of finite character that will be crucial in constructing
our example.\\

\begin{proposition}\label{MR.5} Let $R$ be an integral domain and consider
the map $* :  \bar{F}(R)\rightarrow \bar{F}(R), E\mapsto
E^{*}=\bigcup\{(EJ:J)| J$ runs over the set of finitely generated
fractional ideals of $R$\}. Then $*$ is a semistar operation of
finite character on $R$. Moreover, if $R$ is integrally closed then
$*|_{F(R)}$, is a star operation of finite character on $R$ and
$*$-$Max(R)=Max(R)$.
\end{proposition}

\begin{proof} It is easy to see that $*=d_{a}$, where $d_{a}$ is
the $e. a. b.$ (endlich arithmetisch brauchbar) semistar operation
associated to the trivial semistar operation $d$, \cite[Definition
4.4]{FL2} or \cite[page 4793]{FL}.\\

To see that $*$-$Max(R)=Max(R)$, let $Q$ be a maximal ideal of $R$.
If $Q^{*}=R$, then $1\in Q^{*}$. So there exists a finitely
generated ideal $J$ of $R$ such that $J\subseteq JQ$. Hence $J=JQ$
which contradicts Nakayama's lemma (\cite[Theorem 76]{Ka}). By maximality
$Q=Q^{*}$. Hence $Q$ is $*$-maximal. Conversely, if $Q$ is
$*$-maximal, then $Q$ is prime. So $Q\subseteq M$ for some maximal
ideal $M$ of $R$. By the first step $M$ is $*$-maximal. Hence $Q=M$
and therefore $*$-$Max(R)=Max(R)$.
\end{proof}

\begin{proposition}\label{MR.6} Let $R$ be an integrally closed domain.
Denote by $*_{X}$ the star operation defined on $R[X]$ by
$E^{*_{X}}=\bigcup\{(EJ:J)| J$ runs over the set of finitely
generated fractional ideals of $R[X]$\}, $*_{R}$ the star operation
defined on $R$ by $E^{*_{R}}=\bigcup\{(EJ:J)| J$ runs over the set
of finitely generated fractional ideals of $R$\}. Then:\\
$(1)$ for each nonzero fractional ideal $I$ of $R$,
$(I[X])^{*_{X}}=I^{*_{R}}[X]$.\\
$(2)$ $\overline{*_{X}}=*_{R}$.
\end{proposition}

\begin{proof} Let $I$ be a nonzero fractional ideal of $R$ and
let $f\in I^{*_{R}}[X]$. Write
$f=\displaystyle\sum_{i=0}^{i=n}a_{i}X^{i}$ where $a_{i}\in
I^{*_{R}}$ for each $i$. Then for each $i\in \{0, \dots, n\}$, there
exists a finitely generated ideal $A_{i}$ of $R$ such that
$a_{i}A_{i}\subseteq IA_{i}$. Set
$A=\displaystyle\prod_{i=0}^{i=n}A_{i}$. Clearly $A$ is a finitely
generated ideal of $R$ and $a_{i}A\subseteq IA$ for each $i$. Hence
$fA[X]\subseteq I[X]A[X]$ and so $f\in I[X]^{*_{X}}$.  Then
$I^{*_{R}}[X]\subseteq (I[X])^{*_{X}}$. Conversely, let $f\in
(I[X])^{*_{X}}$. Then there exists a finitely generated ideal $J$ of
$R[X]$ such that $fJ\subseteq I[X]J$. Clearly $c(J)$ is a finitely
generated ideal of $R$ and $c(f)c(J)\subseteq Ic(J)$. Hence
$c(f)\subseteq I^{*_{R}}$. So $f\in I^{*_{R}}[X]$ and therefore
$(I[X])^{*_{X}}=I^{*_{R}}[X]$. Now, for each nonzero fractional
ideal $I$ of $R$, $I^{\overline{*_{X}}}:=(I[X])^{*_{X}}\cap
K=I^{*_{R}}[X]\cap K=I^{*_{R}}$. Hence $\overline{*_{X}}=*_{R}$, as
desired.
\end{proof}

\begin{example}\label{MR.7} The following example shows that:\\
$(1)$ A $*$-maximal ideal $Q$ of $R[X]$ with $Q\cap R=q\not =0$ is
not necessarily extended from $R$, that is, $q[X]\subsetneq Q$.\\
$(2)$ The statements of Theorem~\ref{MR.4} are not equivalent in general.\\
$(3)$ The fact that every upper to zero in $R[X]$ is $*$-invertible
is not enough to ensure that every upper to zero is $*$-maximal
(\cite[Theorem 1.1]{FGH}).\\

Let $k$ be a field and $Y$ an indeterminate over $k$. Let $R=k[Y]$
and let $*_{X}$ be the star operation on $R[X]$ defined in
Proposition~\ref{MR.5}. Let $Q=(X, Y)$. Then $Q\cap R=YR=q$. Since
$*_{X}$-$Max(R[X])= Max(R[X])$, then $Q$ is $*_{X}$-maximal.
However, $q[X]\subsetneq Q$.\\

$(2)$ Since $R$ is a $PID$, then each nonzero ideal of $R$ is divisorial and so $R$ has
exactly one star operation which is the trivial operation $d$. So
$\overline{*_{X}}=d$ and $\overline{*_{X}}$-$Max(R)=Max(R)$. By
\cite[Theorem 19.15]{G1}, for every upper to zero $Q$ in $R[X]$ and
every maximal ideal $M$ of $R$, $Q\not\subseteq M[X]$. However, no
upper to zero $Q$ is $*_{X}$-maximal since $htQ=1$ and all maximal
ideals of $R[X]$ are of height two.\\

$(3)$ Let $Q$ be an upper to zero in $R[X]=k[Y, X]$. Since $R[X]$ is a $UFD$ and $ht(Q)=1$, then $Q$ is principal. Hence
$Q$ is $*_{X}$-invertible but not
$*_{X}$-maximal.
\end{example}

The next results shed light on the relationship between
$*$-$Max(R[X])$, $t$-$Max(R[X])$, $\bar{*}$-$Max(R)$, and
$t$-$Max(R)$ and characterize domains for which every upper to zero
in $R[X]$ is $*$-maximal for a given star operation of finite
character $*$ on $R[X]$.\\

\begin{proposition}\label{MR.8} Let $*$ be a star operation of
finite character on $R[X]$. The following statements are
equivalent:\\
$(i)$ $*$-$Max(R[X])$=$t$-$Max(R[X])$;\\
$(ii)$ $Q_{t}\subsetneq R[X]$ for each $Q\in *$-$Max(R[X])$ such that $Q\cap R\not =0$;\\
$(iii)$ $\bar{*}$-$Max(R)$=$t$-$Max(R)$ and $Q=(Q\cap R)[X]$ for
each $Q\in *$-$Max(R[X])$ with $Q\cap R\not =0$.
\end{proposition}

\begin{proof}
$(ii)\Longrightarrow (i)$ Let $Q\in *$-$Max(R[X])$. If $Q\cap R=0$,
i.e. $Q$ is an upper to zero, then $Q$ is a $t$-prime ideal. So $Q\subseteq
N$ for some $N\in t$-$Max(R[X])$. Since $N$ is a $*$-ideal and $Q$
is $*$-maximal, then $Q=N$. If $Q\cap R\not =0$, by $(ii)$,
$Q_{t}\subsetneq R[X]$. So $Q\subseteq Q_{t}\subseteq N$ for some
$N\in t$-$Max(R[X])$. Since $N$ is a $*$-ideal and $Q$ is
$*$-maximal, then $Q=N$. Hence $*$-$Max(R[X])\subseteq
t$-$Max(R[X])$. Conversely, let $Q\in t$-$Max(R[X])$. Since $Q$ is a
$*$-ideal, then $Q\subseteq N$ for some $N\in *$-$Max(R[X])$. If
$N\cap R=0$, then $N$ is a $t$-ideal and since $Q$ is $t$-maximal,
then $Q=N$. If $N\cap R\not =0$, by $(ii)$, $N_{t}\subsetneq R[X]$.
Since $Q$ is $t$-maximal, then $Q=N$. Hence $Q\in *$-$Max(R[X])$. So
$t$-$Max(R[X])\subseteq
*$-$Max(R[X])$ and therefore $*$-$Max(R[X])$= $t$-$Max(R[X])$.\\

$(i)\Longrightarrow (iii)$ By \cite[Proposition 1.1]{HZ}, $Q=(Q\cap
R)[X]$ for each $Q\in *$-$Max(R[X])$=\\
$t$-$Max(R[X])$ with $Q\cap
R\not =0$.  Let $p\in \bar{*}$-$Max(R)$. Then $(p[X])^{*}\subsetneq
R[X]$. So there exists $Q\in *$-$Max(R[X])$ such that $p[X]\subseteq
(p[X])^{*}\subseteq Q$. By $(i)$, $Q$ is $t$-maximal. Since
$p\subseteq Q\cap R=q$, then $q\in t$-$Max(R)$ and $Q=q[X]$. Since
$q$ is a $\bar{*}$-ideal and $p$ is $\bar{*}$-maximal, then $p=q$.
Hence $p\in t$-$Max(R)$. Conversely, let $p\in t$-$Max(R)$. Since
$p$ is a $\bar{*}$-ideal, then $p\subseteq q$ for some $q\in
\bar{*}$-$Max(R)$. Since $(q[X])^{*}\subsetneq R[X]$, then
$(q[X])^{*}\subseteq Q$ for some $Q\in
*$-$Max(R[X])$=$t$-$Max(R[X])$. Since $p[X]$ is $t$-maximal and
$p[X]\subseteq q[X]\subseteq (q[X])^{*}\subseteq Q$, then $p[X]=Q$.
Hence $p=q\in \bar{*}$-$Max(R)$ and therefore
$\bar{*}$-$Max(R)$=$t$-$Max(R)$.\\

$(iii)\Longrightarrow (ii)$ Let $Q\in *$-$Max(R[X])$ with $Q\cap
R=q\not =0$. By $(ii)$, $Q=q[X]$ and so $q$ is $*$-maximal. Hence
$q$ is $t$-maximal and therefore $Q$ is $t$-maximal. So
$Q_{t}\subsetneq R[X]$, as desired.
\end{proof}

\begin{thm}\label{MR.9} Let $*$ be a star operation of finite character on $R[X]$.
The following statements are equivalent:\\
$i)$ Every upper to zero in $R[X]$ is $*$-maximal;\\
$ii)$ $R$ is $UMT$-domain and $t$-$Max(R[X])$=$*$-$Max(R[X])$.
\end{thm}

\begin{proof} $i)\Longrightarrow ii)$ To see that $R$ is a $UMT$-domain,
let $Q$ be an upper to zero in $R[X]$. Then $Q$ is a
$t$-ideal. So $Q\subseteq N$ for some $t$-maximal ideal $N$ of
$R[X]$. Since $N$ is a $*$-ideal and $Q$ is $*$-maximal, then $Q=N$.
Hence $Q$ is $t$-maximal.\\
Now, let $Q\in *$-$Max(R[X])$ such that $Q\cap R=q\not =(0)$. Then
$Q_{t}\subsetneq R[X]$. For if $Q_{t}=R[X]$, then there exists $g\in
Q$ such that $c(g)_{t}=R$. Let $P$ be a prime ideal of $R[X]$ such
that $gR[X]\subseteq P\subseteq Q$ and $P$ minimal over $gR[X]$.
Then $P$ is $t$-prime (since $P$ is minimal over a principal ideal).
Hence $P\subsetneq Q$. If $P\cap R=0$, then $P$ is an
upper to zero. By $i)$, $P$ is $*$-maximal and so $P=Q$, which is
absurd. Hence $P\cap R=p\not=0$. Let $N$ be a $t$-maximal ideal of
$R[X]$ such that $P\subseteq N$. Since $p\subseteq N\cap R=M$, then
$M$ is a $t$-maximal ideal of $R$ and $N=M[X]$. Since $g\in
P\subseteq N=M[X]$, then $c(g)\subseteq M$. So $R=c(g)_{t}\subseteq
M$, which is absurd. It follows that $Q_{t}\subsetneq R[X]$. By
Proposition~\ref{MR.8},  $t$-$Max(R[X])$=$*$-$Max(R[X])$ as desired.\\
$ii)\Longrightarrow i)$ Trivial.
\end{proof}

We close this section with the following remark-question:\\

\begin{remark} Let $*$ be a star operation of finite character on
$R[X]$ such that every upper to zero is $*$-maximal. By
Theorem~\ref{MR.9} and \cite[Corollary 2.13]{ACo}, it is easy to see that
$w\leq *\leq t$ (here $w$ and $t$ are the $w$- and $t$-operation on $R[X]$).
We are not able to prove or disprove if
the inequalities are strict (i. e., is there a domain $R$ with a star
operation of finite character $*$ on $R[X]$ such that every upper to zero in $R[X]$ is $*$-maximal and $w<*<t$?)
\end{remark}
%%%%%%%%%%%%%%%%%%%%%%%%%%%%%%%%%%%%%%%%%%%%%%%%%%%%%%%%%%%%%%%%%%%%%%%%%%
%%%%%%%%%%%%%%%%%%%%%%%%%%%%%%%%%%%%%%%%%%%%%%%%%%%%%%%%%%%%%%%%%%%%%%%%%%
%%%%%%%%%%%%%%%%%%%%%%%%%%%%%%%%%%%%%%%%%%%%%%%%%%%%%%%%%%%%%%%%%%%%%%%%%%
\section{A Pr\"ufer domain $R$ for which the $v$-operation is not stable}\label{CE}

We start this short section by recalling that a (semi)star operation
on $R$ is stable (respectively finite-stable) if $(A\cap
B)^{*}=A^{*}\cap B^{*}$ for each $A, B$ in $\bar{F}(R)$
(respectively in $f(R)$). More about star operations that distribute
over finite intersection can be found in \cite{AC} and \cite{AC1}.
Our purpose is to answer negatively the following questions
cited in \cite{AC}:\\

\noindent{\bf Question 3.1} (\cite[Question 3.1.(1)]{AC}) Let $*$ be a star operation on a domain $R$.\\
Does $*$ finite-stable $\Longrightarrow$ $*$ is stable?\\

\noindent{\bf Question 3.2} (\cite[Question 3.2.(2)]{AC}) Let $R$ be an integrally closed
domain.\\
Does $v$ finite stable $\Longrightarrow$ $v$ is stable?\\

\noindent{\bf Question 3.3} (\cite[Question 3.3.(1)]{AC}) Is the $v$-operation stable on a
Pr\"ufer domain?\\

{\bf Claim} The $v$-operation on a Pr\"ufer domain $R$ is
finite-stable. Indeed, if $A$ and $B$ are finitely generated ideals
of $R$, then $A$ and $B$ are invertible and so divisorial ideals.
Then $A_{v}\cap B_{v}=A\cap B\subseteq (A\cap B)_{v}\subseteq
A_{v}\cap B_{v}$. Hence $A_{v}\cap B_{v}=A\cap B=(A\cap B)_{v}$.\\

\noindent{\bf Example 3.4} Let $k$ be a field and $X$, $Y$ and $Z$ indeterminates
over $k$. Set $D_{1}=k(X)[Y]_{(Y)}$, $D_{2}=k(Y)[X]_{(X)}$ and
$D=D_{1}\cap D_{2}$. Since $D_{1}$ and $D_{2}$ are independent
valuation domains of $k(X, Y)$ with maximal ideals respectively
$M_{1}=YD_{1}$ and $M_{2}=XD_{2}$, then $D$ is a quasilocal Pr\"ufer
domain with exactly two maximal ideals $Q_{1}=M_{1}\cap D$ and
$Q_{2}=M_{2}\cap D$. Clearly $(D:D_{1})=(D:D_{2})=0$. Set $V=k(X,
Y)[[Z]]=k(X, Y) +M$, $R=D+M$, $A=Z(D_{1}+M)$  and $B=Z(D_{2}+M)$. Then :\\
$i)$ $R$ is a Pr\"ufer domain \cite[Theorem 2.1]{BG}.\\
$(ii)$ Clearly $A^{-1}=B^{-1}=V$ since $(D:D_{1})=(D:D_{2})=0$. Hence $A_{v}=B_{v}=M$ and so $A_{v}\cap B_{v}=M$.\\
$(iii)$ $A\cap B=ZR$, so $(A\cap B)_{v}=ZR\subsetneq ZV=M=A_{v}\cap
B_{v}$.

\bigskip
{\bf Acknowledgment} This work is supported by KFUPM. The author would like to express his sincere thanks to the referee for his/her suggestions and comments.

%%%%%%%%%%%%%%%%%%%%%%%%%%%%%%%%%%%%%%%%%%%%%%%%%%%%%%%%%
%%%%%%%%%%%%%%%%%%%%%%%%%%%%%%%%%%%%%%%%%%%%%%%%%%%%%%%%%
%\bibliographystyle{amsplain}
\bigskip
%%%%%%%%%%%%%%%%%%%%%%%%%%%%%%%%%%%%%%%%%%%%%%%%%%%%%%%%%
%%%%%%%%%%%%%%%%%%%%%%%%%%%%%%%%%%%%%%%%%%%%%%%%%%%%%%%%%
\end{document}